\documentclass[a4paper]{article}
\usepackage{amsfonts}
\usepackage{amsmath}
\usepackage{dsfont}
\usepackage{latexsym}
\usepackage{graphics}
\input epsf

\newcommand\theReals{\mathbb{R}}

\newcommand\theNaturals{\mathbb{N}}

\newcommand\after{\circ}
\newcommand\curly{\{}
\newcommand\ylruc{\}}

\newcommand\intersection{\cap}

\newcommand\interior{\,\mbox{int}}
\newcommand\closure{\,\mbox{cl}}

\newcommand\dee{\,d}

\newtheorem{theorem}{Theorem}
\newtheorem{lemma}{Lemma}

\newtheorem{corollary}{Corollary}

\newtheorem{proposition}{Proposition}

\providecommand\qed{\hfill$\Box$\vspace{10pt}\par}
\renewcommand\qed{\hfill$\Box$\vspace{10pt}\par}

\renewcommand\closure{\mbox{\rm cl }}

\newcommand\poynta{x}
\newcommand\poyntb{y}

\newcommand\topfn{M}

\newcommand\hilbertdist{d_H}
\newcommand\thompsondist{{d_T}}

\title{A Metric Inequality for the Thompson and Hilbert Geometries}
\author{
\parbox{6cm}{
Roger D. {Nussbaum} \\
Mathematics Department \\
Rutgers University \\
New Brunswick \\
NJ 08903. \\
\tt{nussbaum@math.rutgers.edu}}
\parbox{5cm}{
Cormac Walsh \\
INRIA Rocquencourt \\
B.P. 105, \\
78153 Le Chesnay Cedex\\
France. \\
\tt{cormac.walsh@inria.fr}}}
\date{31 March 2004}
 
\begin{document}
 
\maketitle



\begin{abstract}
There are two natural metrics defined on an arbitrary convex cone:
Thompson's part metric and Hilbert's projective metric.
For both, we establish an inequality
giving information about how far the metric is from being
non-positively curved.
\end{abstract}

\section{Introduction}

\newcommand\scalar{\lambda}
\renewcommand\interior{\mbox{int\,}}
\newcommand\oneton{\curly1,\ldots,N\ylruc}

\newcommand\topp{b}
\newcommand\bott{a}
\newcommand\esssup{\mbox{\rm ess sup}}

\newcommand\unit{\mathds{1}}
\newcommand\rel{\equiv}

\newcommand\thompsonnormball{\tilde B^T}
\newcommand\hilbertnormball{\tilde B^H}

Let $C$ be a cone in a vector space $V$.
Then $C$ induces a partial ordering on $V$ given by
$\poynta\le\poyntb$ if and only if $\poyntb-\poynta\in C$.
For each $\poynta\in C\backslash\curly0\ylruc$, $\poyntb\in V$, define
$\topfn(\poyntb/\poynta):=
   \inf\curly\scalar\in\theReals:\poyntb\le\scalar\poynta\ylruc$.
\emph{Thompson's part metric} on $C$ is defined to be
\[
\thompsondist(\poynta,\poyntb)
   := \log\max\Big(\topfn(\poynta/\poyntb),\topfn(\poyntb/\poynta)\Big)
\]
and \emph{Hilbert's projective metric} on $C$ is defined to be
\[
\hilbertdist(\poynta,\poyntb)
   := \log\Big({\topfn(\poynta/\poyntb)}{\topfn(\poyntb/\poynta)}\Big).
\]
Two points in $C$ are said to be in the same part if the distance between
them is finite in the Thompson metric.
If $C$ is~\emph{almost Archimedean}, then,
with respect to this metric, each part of $C$ is
a complete metric space.
Hilbert's projective metric, however, is only a pseudo-metric:
it is possible to find two distinct points which are zero distance apart.
Indeed it is not difficult to see that $\hilbertdist(\poynta,\poyntb)=0$
if and only if $\poynta=\lambda\poyntb$ for some $\lambda>0$.
Thus $\hilbertdist$ is a metric on the space of rays of the cone.
For further details, see Chapter 1 of the monograph~\cite{nussbaum:hilbert}.

Suppose $C$ is finite dimensional and let $S$ be a cross section of $C$,
that is $S:=\curly x\in C:l(x)=1 \ylruc$, where $l:V\to \theReals$ is some
positive linear functional with respect to the ordering on $V$.
Suppose $x,y\in S$ are distinct.
Let $a$ and $b$ be the points in the boundary of $S$
such that $a$, $x$, $y$, and $b$ are collinear and are arranged in this
order along the line in which they lie.
It can be shown that the Hilbert distance between $x$ and $y$ is
then given by the logarithm of the cross ratio of these four points:
\[
\hilbertdist(x,y)=
   \log \frac{|bx|\,|ay|}{|by|\,|ax|}.
\]
Indeed, this was the original definition of  Hilbert.
If $S$ is the open unit disk, the Hilbert metric is exactly the Klein model
of the hyperbolic plane.

An interesting feature of the two metrics above is that they
show many signs of being non-positively curved.
For example, when endowed with the Hilbert metric, the Lorentz cone
$\curly (t,x_1,\dots,x_n)\in\theReals^{n+1}:t^2>x_1^2+\cdots+x_n^2\ylruc$
is isometric to $n$-dimensional hyperbolic space.
At the other extreme, the positive cone
$\theReals_+^n:=
   \curly (x_1,\ldots,x_n):\mbox{$x_i\ge 0$ for $1\le i\le n$}\ylruc$
with either the Thompson or the Hilbert metric is isometric to a
normed space~\cite{harpe}, which one may think of as being flat.
In between, for Hilbert geometries having a strictly-convex $C^2$
boundary with non-vanishing Hessian, the methods of
Finsler geometry~\cite{shenfinsler} apply.
It is known
that such geometries have constant flag curvature $-1$.
More general Hilbert geometries were investigated in~\cite{kellystraus}
where a definition was given of a point of positive curvature.
It was shown that no Hilbert geometries have such points.

However, there are some notions of non-positive curvature which
do not apply. For example,
a Hilbert geometry will only be a CAT(0)
space (see~\cite{bridsonhaefliger})
if the cone is Lorentzian.
Another notion related to negative curvature is that of 
Gromov hyperbolicity~\cite{gromovhyperbolic}.
In~\cite{benoist}, a condition is given characterising those
Hilbert geometries that are Gromov hyperbolic.
This notion has also been investigated in the wider context of
uniform Finsler Hadamard manifolds, which includes certain
Hilbert geometries~\cite{egloff:uniform}.

Busemann has defined non-positive curvature
for~\emph{chord spaces}~\cite{busemann:distinguished}.
These are metric spaces in which there is a distinguished set of
geodesics, satisfying certain axioms.
In such a space, denote by $m_{xy}$ the midpoint along the
distinguished geodesic connecting the pair of points $x$ and $y$.
Then the chord space is non-positively curved if,
for all points $u$, $x$, and $y$ in the space,
\begin{equation}
\label{eqn:riemanninequal}
d(m_{ux},m_{uy})\le \frac{1}{2} d(x,y),
\end{equation}
where $d$ is the metric.

In the case of the Hilbert and Thompson geometries on a part of a closed cone
$C$, there will not necessarily be a unique minimal geodesic connecting
each pair of points.
However, it is known that,
setting $\beta:=M(y/x;C)$ and $\alpha:=1/M(x/y;C)$,
the curve $\phi:[0,1]\to C:$
\begin{equation}
\label{eqn:geodesics}
\phi(s;x,y):=
   \begin{cases}
\displaystyle
   \left(\frac{\beta^s-\alpha^s}{\beta-\alpha}\right)y
      + \left(\frac{\beta\alpha^s-\alpha\beta^s}{\beta-\alpha}\right)x,
           & \mbox{if $\beta\neq\alpha$}, \\
   \alpha^s x,  & \mbox{if $\beta=\alpha$}
   \end{cases}
\end{equation}
is always a minimal geodesic from $x$ to $y$ with respect to both
the Thompson and Hilbert metrics.
We view these as distinguished geodesics.
If the cone $C$ is finite dimensional, then each part of $C$ will be a
chord space under both the Thompson and Hilbert metrics.
Notice that the geodesics above are projective straight lines.
If the cone is strictly convex, these are the only
geodesics that are minimal with respect to the Hilbert metric.
For Thompson's metric, if two points are in the same part of $C$ and are
linearly independent, then there are infinitely many minimal geodesics
between them.

In this paper we investigate whether inequalities
similar to~(\ref{eqn:riemanninequal}) hold for the Hilbert and Thompson
geometries with the geodesics given in~(\ref{eqn:geodesics}).
We prove the following two theorems.

\begin{theorem}
\label{thm:mainthompson}
Let $C$ be an almost Archimedean cone.
Suppose $u,x,y\in C$ are in the same part.
Also suppose that $0<s<1$ and $R>0$, and that
$d_H(u,x)\le R$ and $d_H(u,y)\le R$.
If the linear span of $\curly u,x,y\ylruc$ is $1$- or $2$-dimensional, then
$\thompsondist\big(\phi(s;u,x),\phi(s;u,y))\le s \thompsondist(x,y)$.
In general
\begin{equation}
\label{eqn:thompsonbound}
\thompsondist\big(\phi(s;u,x),\phi(s;u,y)\big)\le
   \left[\frac{2(1-e^{-Rs})}{1-e^{-R}}-s\right]\thompsondist(x,y).
\end{equation}
\end{theorem}
Note that the bracketed value on the right hand
side of this inequality is strictly increasing in $R$.
As $R\to 0$, this value goes to $s$, which reflects the fact that
in small neighborhoods the Thompson metric looks like a norm.
As $R\to\infty$, the bracketed value goes to $2-s$.

\begin{theorem}
\label{thm:mainhilbert}
Let $C$ be an almost Archimedean cone.
Suppose $u,x,y\in C$ are in the same part.
Also suppose that $0<s<1$ and $R>0$ and that
$d_H(u,x)\le R$ and $d_H(u,y)\le R$.
If the linear span of $\curly u,x,y\ylruc$ is $1$- or $2$-dimensional, then
$\hilbertdist\big(\phi(s;u,x),\phi(s;u,y))\le s \hilbertdist(x,y)$.
In general
\begin{equation}
\label{eqn:hilbertbound}
\hilbertdist\big(\phi(s;u,x),\phi(s;u,y)\big)\le
   \left[\frac{1-e^{-Rs}}{1-e^{-R}}\right]\hilbertdist(x,y).
\end{equation}
\end{theorem}
Again, the bracketed value on the right hand side
increases strictly with increasing $R$.
This time, it ranges between $s$ as $R\to 0$ and
$1$ as $R\to\infty$.

Our method of proof will be to first establish the results when
$C$ is the positive cone $\theReals_+^N$, with $N\ge 3$.
It will be obvious from the proofs that the bounds given are the
best possible in this case.
A crucial lemma will state that any finite set of $n$ elements
of a Thomson or Hilbert geometry can be isometrically embedded in
$\theReals_+^{n(n-1)}$ with, respectively, its Thompson
or Hilbert metric.
This lemma will allow us to extend the same bounds to more general cones,
although in the general case the bounds may no longer be tight.

A special case of Theorem~\ref{thm:mainhilbert}
was proved in~\cite{socie_behaviour} using a simple geometrical argument.
It was shown that if two particles start at the same point and
travel along distinct straight-line geodesics at unit speed in the Hilbert
metric, then the Hilbert distance between them is strictly increasing.
This is equivalent to the special case of Theorem~\ref{thm:mainhilbert}
when $d_H(u,x)=d_H(u,y)$ and $R$ approaches infinity.

A consequence of Theorems~\ref{thm:mainthompson} and \ref{thm:mainhilbert}
is that both the Thompson and Hilbert geometries are semihyperbolic
in the sense of Alonso and Bridson~\cite{semihyperbolic}.
Recall that a metric space is semihyperbolic if it admits a bounded
quasi-geodesic bicombing.
A bicombing is a choice of path between each pair of points.
We may use the one given by
\[
\zeta_{(x,y)}(t):=
\begin{cases}
\displaystyle
\phi\Big(\frac{t}{d(x,y)},x,y\Big), & \mbox{if $t\in[0,d(x,y)]$} \\
y, & \mbox{otherwise}
\end{cases}
\]
for each pair of points $x$ and $y$ in the same part of $C$.
Here $d$ is either the Thompson or Hilbert metric.
This bicombing is geodesic and hence quasi-geodesic.
To say it is bounded means that there exist constants
$M$ and $\epsilon$ such that
\[
d(\zeta_{(x,y)}(t),\zeta_{(w,z)}(t))\le M\max(d(x,w),d(y,z))+\epsilon
\]
for each $x,y,w,z\in C$ and $t\in[0,\infty)$.
\begin{corollary}
\label{cor:thompson}
Each part of $C$ is semihyperbolic when endowed with either Thompson's
part metric or Hilbert's projective metric.
\end{corollary}

It should be pointed out that for some cones there are other good
choices of distinguished geodesics.
For example, for the cone of positive definite symmetric matrices
$\mbox{Sym}(n)$, a natural choice would be 
$\phi(s;X,Y):= X^{1/2}(X^{-1/2}YX^{-1/2})^sX^{1/2}$
for $X,Y\in\mbox{Sym}(n)$ and $s\in[0,1]$.
It can be shown that, with this choice, $\mbox{Sym}(n)$ is non-positively
curved in the sense of Busemann under both the Thompson and Hilbert
metrics.
This result has been generalized to both
symmetric cones~\cite{walshgunawardena}
and to the cone of positive elements of a $C^*$-algebra~\cite{cpr:convexity}.

Although Hilbert's projective metric arose in geometry, it has also
been of great interest to analysts. This is because many naturally
occurring maps in analysis, both linear and non-linear, are
either non-expansive or contractive with respect to it. Perhaps the first
example of this is due to
G.~Birkhoff~\cite{birkhoff:jentzsch,birkhoff:uniformly},
who noted that matrices with strictly positive entries (or indeed
integral operators with strictly positive kernels)
are strict contractions with respect to Hilbert's metric.
References to the literature connecting this metric to
positive linear operators can be found
in~\cite{EvesonNussbaum:elementary,EvesonNussbaum:Applications}.
It has also been used to study the spectral radii of elements
of Coxeter groups~\cite{mcmullen}.
Both metrics have been applied to questions concerning the convergence
of iterates of non-linear
operators~\cite{bushell:hilbert,walshgunawardena,nussbaum:hilbert,
nussbaum:iterated2,nussbaum:omega}.
The two metrics have been used to solve problems involving
non-linear integral equations~\cite{potter:applications,thompson:certain},
linear operator equations~\cite{bushell:hilbert,bushell:cayley}, and
ordinary differential equations~\cite{birkhoff:integro,nussbaum:omega,
wysocki:behavior,wysocki:ergodic}.
Thompson's metric has also been usefully applied
in~\cite{nussbaum:iterated2, nussbaum:entropy} to obtain
``DAD theorems'', which are scaling results concerning kernels of
integral operators. Another application of this metric is in
Optimal Filtering~\cite{liveraniwojtkowski},
while Hilbert's metric has been used in
Ergodic Theory~\cite{liverani} and Fractal Diffusions~\cite{metz}.

\section{Proofs}

A cone is a subset of a (real) vector space that is convex,
closed under multiplication by positive scalars,
and does not contain any vector subspaces of dimension one.
We say that a cone is almost Archimedean if the closure
of its restriction to any two-dimensional subspace is also a cone.

The proofs of Theorems~\ref{thm:mainthompson} and ~\ref{thm:mainhilbert}
will involve the use of some infinitesimal arguments.
We recall that both the Thompson and Hilbert geometries are
\emph{Finsler} spaces~\cite{nussbaum:finsler}.
If $C$ is a closed cone in $\theReals^N$ with non-empty interior, then
$\interior C$ can be considered to be an $N$-dimensional manifold
and its tangent space at each point can be identified with $\theReals^N$.
If a norm
\[
|v|^T_x := \inf\curly\alpha>0:-\alpha x\le v\le \alpha x\ylruc
\]
is defined on the tangent space at each point $x\in \interior C$,
then the length of any piecewise $C^1$ curve $\alpha:[a,b]\to\interior C$
can be defined to be
\[
L^T(\alpha):=\int_a^b|\alpha'(t)|^T_{\alpha(t)}\dee t.
\]
The Thompson distance between any two points is recovered by
minimizing over all paths connecting the points:
\[
\thompsondist(x,y)=
   \inf\curly
      L^T(\alpha): \alpha\in PC^1[x,y]
   \ylruc,
\]
where $PC^1[x,y]$ denotes the set of all piecewise $C^1$ paths
$\alpha:[0,1]\to\interior C$ with $\alpha(0)=x$ and $\alpha(1)=y$.
A similar procedure yields the Hilbert metric when the norm above
is replaced by the semi-norm
\[
|v|^H_x := M(v/x) - m(v/x).
\]
Here $M(v/x)$ is as before and
$m(v/x):= \sup\curly\scalar\in\theReals:v\ge\scalar x\ylruc$.
The Hilbert geometry will be Riemannian only in the case of the
Lorentz cone.  The Thompson geometry will be Riemannian
only in the trivial case of the one-dimensional cone $\theReals_+$.

Our strategy will be to first prove the theorems for the case of the
positive cone $\theReals_+^N$, and then extend them to the general case.
The proof in the case of $\theReals_+^N$ will involve investigation
of the map $g:\interior\theReals_+^N\to\interior\theReals_+^N$:
\begin{equation}
\label{eqn:geedef}
g(x):=\phi(s;\unit,x)=
   \begin{cases}
\displaystyle
   \left(\frac{\topp^s-\bott^s}{\topp-\bott}\right)x
      + \left(\frac{\topp\bott^s-\bott\topp^s}{\topp-\bott}\right)\unit,
           & \mbox{if $\topp\neq\bott$}, \\
   \bott^s \unit,  & \mbox{if $\topp=\bott$},
   \end{cases}
\end{equation}
where $\topp:=\topp(x):=\max_i x_i$ and $\bott:=\bott(x):=\min_i x_i$.
Here $s\in(0,1)$ is fixed and we are using the notation
$\unit:=(1,\ldots,1)$.
The derivative of $g$ at $x\in\interior\theReals_+^N$ is
a linear map from $\theReals^N\to\theReals^N$. Taking $|\cdot|_x^T$ as
norm on the domain and $|\cdot|_{g(x)}^T$ as norm on the range,
the norm of $g'(x)$ is
\[
||g'(x)||_T := \sup\curly|g'(x)(v)|_{g(x)}^T:|v|_x^T\le 1\ylruc.
\]
If, instead, we take the appropriate infinitesimal Hilbert semi-norms
on the domain and range, then the norm of $g'(x)$ is given by
\[
||g'(x)||_H := \sup\curly|g'(x)(v)|_{g(x)}^H:|v|_x^H\le 1\ylruc.
\]

For each pair of distinct integers $I$ and $J$ contained in $\{1,\ldots,N\}$,
let
\[
U_{I,J}:=
   \Big\{x\in\interior\theReals_+^N:
      \mbox{$0<x_I<x_i<x_J$ for all $i\in\{1,\ldots,N\}\backslash\{I,J\}$}
   \Big\}.
\]
On each set $U_{I,J}$, the map $g$ is $C^1$ and is given by the formula
\[
g(x)=   \left(\frac{x_J^s-x_I^s}{x_J-x_I}\right)x
      + \left(\frac{x_J x_I^s-x_Ix_J^s}{x_J-x_I}\right)\unit.
\]
Let $U$ denote the union of the sets $U_{I,J}$; $I,J\in\{1,\ldots,N\}$,
$I\neq J$. If $x\in\theReals_+^N\backslash U$, then there must exist distinct
integers $m,n\in\{1,\ldots,N\}$ with either $x_n=x_m=\max_ix_i$
or $x_n=x_m=\min_ix_i$.
The set $x\in\theReals_+^N$ with $x_n=x_m$ has ($N$-dimensional) Lebesgue
measure zero, so the complement of $U$ in $\theReals_+^N$ has
Lebesgue measure zero.

We recall the following results from~\cite{nussbaum:finsler}.
The first is a combination of Corollaries 1.3 and 1.5 from that paper.
\begin{proposition}
\label{prp:thompson}
Let $C$ be a closed cone with non-empty interior
in a finite dimensional normed space $V$.
Suppose $G$ is an open subset of $\interior C$ such that
$\phi(s;x,y)\in G$ for all $x,y\in G$ and $s\in[0,1]$.
Suppose also that $f:G\to \interior C$ is a locally Lipschitzian map
with respect to the norm on $V$. 
Then
\[
\inf\curly k\ge0:
   \mbox{$\thompsondist(f(x),f(y))\le k\thompsondist(x,y)$ for all $x,y\in G$}
   \ylruc
      = \displaystyle\esssup_{x\in G}
||f'(x)||_T.
\]
\end{proposition}
It is useful in this context to recall that every locally
Lipschitzian map is Fr\'{e}chet differentiable Lebesgue almost everywhere.
The next proposition is a special case of Theorem 2.5
in~\cite{nussbaum:finsler}.
\begin{proposition}
\label{prp:hilbert}
Let $C$ be a closed cone with non-empty interior in a
normed space $V$ of finite dimension $N$.
Let $l$ be a linear functional on $V$ such that $l(x)>0$
for all $x\in \interior C$, and define $S:=\curly x\in C:l(x)=1\ylruc$.
Let $G$ be a relatively-open convex subset of S.
Suppose that $f:G\to \interior C$
is a locally Lipschitzian map with respect to the norm
on $V$.
Then
\[
\inf\curly k\ge0:
   \mbox{$\hilbertdist(f(x),f(y))\le k\hilbertdist(x,y)$ for all $x,y\in G$}
   \ylruc
      = \displaystyle\esssup_{x\in G}
||f'(x)||_{\tilde H},
\]
where
$||f'(x)||_{\tilde H}
   := \sup\curly|f'(x)(v)|_{f(x)}^H:\mbox{$|v|_x^H\le 1$, $l(v)=0$}\ylruc$.
Here the essential supremum is taken with respect to the
$N-1$-dimensional Lebesgue measure on $S$.
\end{proposition}

Since we wish to apply Propositions~\ref{prp:thompson} and~\ref{prp:hilbert}
to the map $g$, we must prove that it is locally Lipschitzian.
\begin{lemma}
\label{lem:lipschitz}
The map $g:\interior(\theReals_+^N)\to\interior(\theReals_+^N)$
defined by~(\ref{eqn:geedef}) is locally Lipschitzian.
\end{lemma}
\begin{proof}
We use the supremum norm $||x||_\infty:=\max_i|x_i|$ on $\theReals^N$.
Clearly, $|b(x)-b(y)|\le||x-y||_\infty$ and
$|a(x)-a(y)|\le||x-y||_\infty$ for all $x,y\in\interior(\theReals_+^N)$.
Therefore both $a$ and $b$ are Lipschitzian with Lipschitz constant 1.

Let $\gamma:[0,\infty)\to[0,\infty)$ be defined by
\[
\gamma(t):=
   \begin{cases}
   \displaystyle\frac{t^s-1}{t-1}, & \mbox{for $t\neq1$}, \\
   s, & \mbox{for $t=1$}.
   \end{cases}
\]
Then $g$ may be expressed as
\[
g(x)=
   a^{s-1}\gamma\left({b/a}\right)x
       + a^s\Big(1-\gamma\left({b/a}\right)\Big)\unit.
\]
The Binomial Theorem gives that
\[
\gamma(t) = \sum_{k=1}^\infty
   \binom{s}{k}(t-1)^k
\qquad\mbox{for $|t-1|<1$}
\]
and so $\gamma$ is $C^\infty$ on a neighborhood of 1.
Hence it is $C^\infty$ on $[0,\infty)$, and thus locally Lipschitzian.
It follows that $g$ is also locally Lipschitzian.
\qed
\end{proof}

\subsection{Thompson's metric}

We have the following bound on the norm of $g'(x)$
with respect to the Thompson metric.
\begin{lemma}
\label{lem:thompsonnorm}
Consider the Thompson metric on $\interior\theReals_+^N$.
Let $x\in U_{1,N}$.
If $N=1$ or $N=2$ then the norm of $g'$ at $x$ is given by
$||g'(x)||_T=s$.
If $N\ge3$, then
\begin{equation}
\label{eqn:normexpression}
||g'(x)||_T  =
        \frac{x_N-x_{N-1}}{x_N-x_1}
        \theta\left(\frac{x_N}{x_1}\right)
        \frac{x_1^{s+1}}{E_{N-1}}
     +\frac{(x_N^s-x_1^s)x_{N-1}}{E_{N-1}}
     + \frac{x_{N-1}-x_1}{x_N-x_1}
        \theta\left(\frac{x_1}{x_N}\right)
        \frac{x_N^{s+1}}{E_{N-1}}
\end{equation}
where $\theta(t):=(1-s)-t^s+st$
and $E_i(x):=E_i:=x_i(x_N^s-x_1^s)+x_Nx_1^s-x_1x_N^s$.
\end{lemma}
\begin{proof}
If $N=1$ and $x>0$, then $g(x)=x^s$. We leave the proof in this case to
the reader and assume that $N\ge2$.

For $x\in U_{1,N}$,
\[
g(x)=   \left(\frac{x_N^s-x_1^s}{x_N-x_1}\right)x
      + \left(\frac{x_N x_1^s-x_1x_N^s}{x_N-x_1}\right)\unit.
\]
Let
\[
h_{ij}(x) := \frac{x_j}{g_i(x)}\frac{\partial g_i}{\partial x_j}(x).
\]
Straightforward calculation gives,
for each $j\in\{1,\ldots,N\}$,
\begin{eqnarray*}
h_{1j}(x) &=& s\delta_{1j} \\
\mbox{and}\qquad
h_{Nj}(x) &=& s\delta_{Nj}.
\end{eqnarray*}
Here $\delta_{ij}$ is the Kronecker delta function which takes the value
$1$ if $i=j$ and the value $0$ if $i\neq j$.
Clearly, $h_{ij}(x)=0$ for $1<i<N$ and $j\not\in\{1,i,N\}$.
For $1<i<N$,
\begin{eqnarray}
\label{eqn:hachesstart}
h_{i1}(x) &=& \frac{x_N-x_i}{x_N-x_1}
                   \, \theta\left(\frac{x_N}{x_1}\right)
              \frac{x_1^{s+1}}{E_i}
    \qquad\ge 0,\\
\label{eqn:hachesmid}
h_{ii}(x) &=& \frac{x_N^s-x_1^s}{E_i}x_i
    \qquad\qquad\qquad\qquad\ge 0,\\
\label{eqn:hachesend}
h_{iN}(x) &=& -\frac{x_i-x_1}{x_N-x_1}
                   \, \theta\left(\frac{x_1}{x_N}\right)
              \frac{x_N^{s+1}}{E_i}
    \qquad\le 0.
\end{eqnarray}
Inequalities~(\ref{eqn:hachesstart}--\ref{eqn:hachesend})
rely on the fact that
$\theta(t)\ge 0$ for $t\ge0$. This may be established by observing that
$\theta(1)=\theta'(1)=0$ and $\theta''(t)>0$ for $t\ge0$.

Let
\[
\thompsonnormball:=\left\{v\in\theReals^N:
           \mbox{$\max_j |v_j| \le 1$}\right\}.
\]
We wish to calculate
\begin{equation}
\label{eqn:gnorm}
||g'(x)||_T = \sup\bigg\{\Big|\sum_j h_{ij}v_j\Big|:
     \mbox{$1\le i\le N$, $v\in \thompsonnormball$}\bigg\}.
\end{equation}

For $i=1$ or $i=N$, we have $|\sum_jh_{ij}v_j|\le s$
for any choice of $v\in \thompsonnormball$.
If $N=2$, then it follows that $||g'(x)||_T=s$ for all $x\in U_{1,N}$.

For the rest of the proof we shall therefore assume that $N\ge3$.
For $1<i<N$, it is clear from
Inequalities~(\ref{eqn:hachesstart}-\ref{eqn:hachesend})
that $|\sum_jh_{ij}v_j|$ is maximized when $v_1=v_i=1$ and $v_N=-1$.
In this case
\begin{eqnarray}
\label{eqn:thompsonnorm}
\Big|\sum_jh_{ij}v_j\Big| &=&
        \frac{1}{E_i} \Big[
        \frac{x_N-x_i}{x_N-x_1}
\theta\left(\frac{x_N}{x_1}\right)x_1^{s+1}
             + (x_N^s-x_1^s)x_i
             + \frac{x_i-x_1}{x_N-x_1}
                  \theta\left(\frac{x_1}{x_N}\right)x_N^{s+1}
          \Big] \\
\label{eqn:ratfunc}
   &=& \frac{c_1x_i+c_2}{c_3x_i+c_4},
\end{eqnarray}
where $c_1$, $c_2$, $c_3$, and $c_4$ depend on $x_1$ and $x_N$
but not on $x_i$.
Observe that $c_3x_i+c_4\neq0$ for $x_1\le x_i\le x_N$.
Given this fact, the general form of Expression~(\ref{eqn:ratfunc})
leads us to conclude that it is either non-increasing or non-decreasing
when regarded as a function of $x_i$.
When we substitute $x_i=x_1$, we get $|\sum_jh_{ij}v_j|=s$.
When we substitute $x_i=x_N$, we get
\begin{equation}
\Big|\sum_jh_{ij}v_j\Big|=
   \frac{2\Big(1-({x_1}/{x_N})^s\Big)}{1-({x_1}/{x_N})}-s.
\end{equation}
Now, writing $\Gamma(t):=2(1-t^s)/(1-t)-s$, we have
$\Gamma'(t)=-2t^s\theta(t^{-1})/(1-t)^2<0$,
in other words $\Gamma$ is decreasing on $(0,1)$.
In particular, $\Gamma(x_1/x_N)\ge\lim_{t\to1}\Gamma(t)=s$.
Therefore Expression~(\ref{eqn:thompsonnorm})
is non-decreasing in $x_i$.
So, the supremum in~(\ref{eqn:gnorm}) is attained when $v$ is as above
and $i=N-1$. Recall that $x_{N-1}$ is the second largest component of $x$.
The conclusion follows.
\qed
\end{proof}
\begin{corollary}
\label{cor:thompsonnorm}
Let $R>0$. If $N=1$ or $N=2$, then
$\esssup \curly||g'(x)||_T:x\in \interior\theReals_+^N\ylruc=s$.
If $N\ge3$, then
\[
\esssup\curly||g'(x)||_T:d_H(x,\unit)\le R\ylruc 
   = \frac{2(1-e^{-Rs})}{1-e^{-R}}-s.
\]
\end{corollary}
\begin{proof}
Note that if $\sigma:\theReals_+^N\to\theReals_+^N$ is some permutation
of the components, then $g\after\sigma(x)=\sigma\after g(x)$
for all $x\in\theReals_+^N$. Furthermore, $\sigma$ will be
an isometry of both the Thompson and Hilbert metrics.
It follows that, given any $x\in U_{I,J}$ with $I,J\in\{1,\ldots,N\}$,
$I\neq J$, we may reorder the components of $x$ to find a point $y$ in
$U_{1,N}$ such that $||g'(y)||_T=||g'(x)||_T$.
Recall, also, that the complement of $U$ in $\interior\theReals_+^N$
has $N$-dimensional Lebesgue measure zero.
From these two facts, it follows that
the essential supremum of $||g'(x)||_T$ over
$\overline B_R(\unit):=
   \curly x\in \interior\theReals_+^N:d_H(x,\unit)\le R\ylruc$
is the same as its supremum over
$\overline B_R(\unit)\intersection U_{1,N}$.

In the case when $N=1$ or $N=2$, the conclusion follows immediately.

For $N=3$, we must maximize Expression~(\ref{eqn:normexpression}) under the
constraints $x_1<x_{N-1}<x_N$ and $x_1/x_N\ge\exp(-R)$.
First, we maximize over $x_{N-1}$, keeping $x_1$ and $x_{N}$ fixed.
In the proof of the previous lemma,
we showed that Expression~(\ref{eqn:normexpression}) is non-decreasing
in $x_{N-1}$, and so it will be maximized when $x_{N-1}$ approaches
$x_{N}$. Here it will attain the value
\begin{equation}
\label{eqn:thompsongammax1xN}
   \frac{2\Big(1-({x_1}/{x_N})^s\Big)}{1-({x_1}/{x_N})}-s
  = \Gamma(x_1/x_N).
\end{equation}
We also showed that $\Gamma$ is decreasing on $(0,1)$.
Therefore~(\ref{eqn:thompsongammax1xN})
will be maximized when $x_1/x_N=\exp(-R)$,
where it takes the value
\[
   \frac{2(1-e^{-Rs})}{1-e^{-R}}-s.
\]
\qed
\end{proof}

\begin{lemma}
\label{lem:reduce}
Let $C$ be an almost Archimedean cone and let
$\curly x_i:i\in I\ylruc$
be a finite collection of elements of $C$ of cardinality $n$, all lying in
the same part.
Denote by $W$ the linear span of $\curly x_i:i\in I\ylruc$
and write $C_W:=C\intersection W$.
Denote by $\interior C_W$ the interior of $C_W$ as a subset of $W$,
using on $W$ the unique Hausdorff linear topology.
Then each of the points $x_i;i\in I$ is contained in $\interior C_W$.
Furthermore, there exists a linear map $F:W\to\theReals^{n(n-1)}$
such that
$F(\interior C_W)\subset\interior\theReals_+^{n(n-1)}$ and
\begin{equation}
\label{eqn:equalm}
M(x_i/x_j;C)=M(F(x_i)/F(x_j);\theReals_+^{n(n-1)})
\end{equation}
for each $i,j\in I$.
\end{lemma}
\begin{proof}
Since the points $\curly x_i:i\in I\ylruc$ all lie
in the same part of $C$, they also all lie in the same part of $C_W$.
Therefore there exist positive constants $a_{ij}$
such that $x_j-a_{ij}x_i\in C_W$ for all $i,j\in I$.
If we define $a:=\min\curly a_{ij}:i,j\in I\ylruc$
it follows that $x_j+\delta x_i\in C_W$
whenever $|\delta|\le a$ and $i,j\in I$.
Now select $i_1,\ldots,i_m\in I$ such that $\curly x_{i_k}:1\le k\le m\ylruc$
form a linear basis for $W$.
For each $y\in W$, we define $||y||:=\max\curly|b_k|:1\le k\le m\ylruc$,
where $y=\sum_{k=1}^m b_kx_{i_k}$ is the unique representation
of $y$ in terms of this basis.
The topology on $W$ generated by this norm is the same as the
one we have been using.
If $||y||\le a/m$ and $j\in I$, then $x_j+mb_k x_{i_k}\in C_W$
for $1\le k \le m$. It follows that
\[
x_j+y=\frac{1}{m}\sum_{k=1}^m(x_j+mb_kx_{i_k})\in C_W
\]
whenever $||y||\le a/m$. This proves that $x_j\in \interior C_W$
for all $j\in I$.

It is easy to see that $\beta_{ij}:= M(x_i/x_j;C)=M(x_i/x_j;C_W)$
for all $i,j\in I$, $i\neq j$.
Observe that $\beta_{ij}x_j-x_i\in\partial C_W$.
Since $\interior C_W$ is a non-empty open convex set which does not
contain $\beta_{ij}x_j-x_i$,
the geometric version of the Hahn-Banach Theorem implies that there exists a
linear functional $f_{ij}:W\to\theReals$ and a real number
$r_{ij}$ such that $f_{ij}(\beta_{ij}x_j-x_i)\le r_{ij}<f_{ij}(z)$
for all $z\in\interior C_W$. Because $0$ is in the closure of $\interior C_W$
and $f_{ij}(0)=0$, we have $r_{ij}\le0$.
On the other hand, if $f_{ij}(z)<0$
for some $z\in\interior C_W$, then considering $f_{ij}(tz)$ we see that
$f_{ij}$ would not be bounded below on $\interior C_W$. It follows that
$r_{ij}=0$. Since $\beta_{ij}x_j-x_i$ is in the closure of $\interior C_W$,
we must have $f_{ij}(\beta_{ij}x_j-x_i)=0$.

Now, define
\[
F:W\to\theReals^{n(n-1)}:z\mapsto(f_{ij}(z))_{i,j\in I, \, i\neq j},
\]
so that $f_{ij}(z);i,j\in I, i\neq j$ are the components of $F(z)$.
Clearly, $F$ is linear and maps $\interior C_W$
into $\interior\theReals_+^{n(n-1)}$.
Also, for all $i,j\in I$, $i\neq j$,
\[
M(F(x_i)/F(x_j);\theReals_+^{n(n-1)})
   = \inf\curly\lambda>0:
            \mbox{$f_{kl}(\lambda x_j-x_i)\ge 0$ for all $k,l\in I$, $k\neq l$}
       \ylruc.
\]
For $\lambda\ge\beta_{ij}$, we have $\lambda x_j-x_i\in \closure C_W$ and so
$f_{kl}(\lambda x_j-x_i)\ge 0$ for all $k,l\in I$, $k\neq l$.
On the other hand, for $\lambda<\beta_{ij}$,
we have $f_{ij}(\lambda x_j-x_i)<0$
since $f_{ij}(x_j)>0$.
We conclude that $M(F(x_i)/F(x_j);\theReals_+^{n(n-1)})= \beta_{ij}$.
\qed
\end{proof}

\begin{lemma}
\label{lem:thompsonpos}
Theorem~\ref{thm:mainthompson} holds in the special case when
$C=\theReals_+^N$ with $N\ge 3$.
\end{lemma}
\begin{proof}
Each part of $\theReals_+^N$ consists of elements of $\theReals_+^N$
all having the same components equal to zero.
Thus each part can be naturally identified
with $\interior\theReals_+^n$, where $n$ is the number of strictly
positive components of its elements. We may therefore assume initially that
$\curly x,y,u\ylruc\subset \interior\theReals_+^N$.

Define $L:\theReals^N\to\theReals^N$ by
$L(z):=(u_1z_1,\ldots,u_Nz_N)$.
Its inverse is given by $L^{-1}(z):=(u^{-1}_1z_1,\ldots,u^{-1}_Nz_N)$.
Both $L$ and $L^{-1}$ are linear maps which leave $\theReals_+^N$
invariant. It follows that $L$ and $L^{-1}$ are isometries of $\theReals_+^N$
with respect to both the Thompson and Hilbert metrics.
Therefore, for $u,z\in\interior\theReals_+^N$,
\[
L^{-1}(\phi(s;u,z))=\phi(s;L^{-1}(u),L^{-1}(z)).
\]
Thus, we may as well assume that $u=\unit$.

We now wish to apply Proposition~\ref{prp:thompson}
with $f:=g$ and
$G:=B_{R+\epsilon}(\unit)
   =\curly z\in\theReals_+^N:d_H(z,\unit)< R+\epsilon\ylruc$.
It was shown in~\cite{nussbaum:hilbert} that $G$ is a convex cone,
in other words that it is closed under multiplication by positive scalars and
under addition of its elements.
Since $\phi(s;w,z)$ is a positive combination of $w$ and $z$,
it follows that $\phi(s;w,z)$ is in $G$ if $w$ and $z$ are.
If we now apply Lemma~\ref{lem:lipschitz},
Proposition~\ref{prp:thompson},
and Corollary~\ref{cor:thompsonnorm}, and let $\epsilon$ approach zero,
we obtain the desired result.
\qed
\end{proof}

\begin{lemma}
\label{lem:thompson2d}
Theorem~\ref{thm:mainthompson} holds in the special case when
the linear span of $\curly x,y,u\ylruc$ is one- or two-dimensional.
\end{lemma}
\begin{proof}
Let $W$ denote the linear span of $\curly x,y,u\ylruc$, in other words
the smallest linear subspace containing these points.
By Lemma~\ref{lem:reduce},
$x$, $y$, and $u$ are in the interior of $C\intersection W$ in $W$.
It is easy to see that
$M(z/w;C)=M(z/w;C\intersection W)$ for all
$w,z\in \interior (C\intersection W)$.
Therefore, we can work in the cone $C\intersection W$.

It is not difficult to show~\cite{EvesonNussbaum:elementary} that
if $m:=\dim W$ is either one or two,
then there is a linear isomorphism $F$ from $W$ to
$\theReals^m$ taking $\interior (C\intersection W)$
to $\interior\theReals_+^m$.
It follows that $F$ is an isometry of both the Thompson and Hilbert metrics
and $F(\phi(s;z,w))=\phi(s;F(z),F(w))$
for all $z,w\in\interior(C\intersection W)$.
We may thus assume that $C=\theReals_+^m$ and
$u,x,y\in \interior C$.

As in the proof of Lemma~\ref{lem:thompsonpos},
we may assume that $u=\unit$.

To obtain the required result,
we apply Lemma~\ref{lem:lipschitz},
Corollary~\ref{cor:thompsonnorm},
and Proposition~\ref{prp:thompson}
with $f:=g$ and $G:=\interior \theReals_+^m$.
\qed
\end{proof}

\begin{proof}[of Theorem~\ref{thm:mainthompson}]
Let $W$ denote the linear span of
$\curly x,y,u\ylruc$.
Lemma~\ref{lem:thompson2d} handles the case when
these three points are not linearly independent;
we will therefore assume that they are.
Thus the five points $x$, $y$, $u$, $\phi(s;u,x)$, and $\phi(s;u,y)$
are distinct.
We apply Lemma~\ref{lem:reduce} and obtain a linear
map $F:W\to\theReals_+^{20}$ with the specified properties.
From~(\ref{eqn:equalm}), it is clear that
$\thompsondist(z,w)=\thompsondist'(F(z),F(w))$ for each
$z,w\in\curly x,y,u,\phi(s;u,x),\phi(s;u,y)\ylruc$.
Here we are using $\thompsondist'$ to denote the Thompson metric
on $\theReals_+^{20}$.
Note that $\phi(s;u,x)$ is a positive combination of $u$ and $x$ and that
the coefficients of $u$ and $x$ depend only on $s$, $M(u/x;C)$,
and $M(x/u;C)$.
The latter two quantities are equal to $M(F(u)/F(x);\theReals_+^{20})$
and $M(F(x)/F(u);\theReals_+^{20})$ respectively.
We conclude that $F(\phi(s;u,x))=\phi(s;F(u),F(x))$.
A similar argument gives $F(\phi(s;u,y))=\phi(s;F(u),F(y))$.
Inequality~(\ref{eqn:thompsonbound}) follows by applying
Lemma~\ref{lem:thompsonpos} to the points $F(x)$, $F(y)$, and $F(u)$
in the cone $\theReals_+^{20}$.
\qed
\end{proof}

\subsection{Hilbert's metric}

We shall continue to use the same notation.
Thus, for a given $N\in\theNaturals$ and
$s\in(0,1)$, we use $g$ to denote the function in~(\ref{eqn:geedef})
and $U$ to denote the union of sets $U_{I,J}$
with $I,J\in\curly1,\ldots,N\ylruc$, $I\neq J$.
We also use the functions $\theta(t):=(1-s)-t^s+st$
and $E_i(x):=E_i:=x_i(x_N^s-x_1^s)+x_Nx_1^s-x_1x_N^s$, and write
$h_{ij}(x) := ({x_j}/{g_i(x)}){\partial g_i}/{\partial x_j}(x)$.
As was noted earlier, $\theta(t)>0$ if $t>0$ and $t\neq1$. Also,
$\gamma(t):= (1-t^s)/(1-t)$, $\gamma(1):=s$ is strictly decreasing on
$[0,\infty)$.
We shall also use the simple but useful observation that
if $c_1$, $c_2$, $c_3$, and $c_4$ are constants such that
$c_3t+c_4\neq0$ for $a\le t\le b$, then the function
$t\mapsto (c_1t+c_2)/(c_3t+c_4)$
is either increasing on $[a,b]$ (if $c_1c_4-c_2c_3\ge0$)
or decreasing on $[a,b]$ (if $c_1c_4-c_2c_3\le0$).
Either way, the function attains is maximum over $[a,b]$
at $a$ or $b$.

Recall that if $g$ is Fr\'{e}chet differentiable at
$x\in\interior\theReals_+^N$
then $||g'(x)||_H$ denotes the norm of $g'(x)$ as a linear map
from $(\theReals^N,||\cdot||_x^H)$ to 
$(\theReals^N,||\cdot||_{g(x)}^H)$, although, of course, $||\cdot||_x^H$
and $||\cdot||_{g(x)}^H$ are semi-norms rather than norms.

\begin{lemma}
Consider the Hilbert metric on $\interior\theReals_+^N$ with $N\ge 2$.
Let $x\in U_{1,N}$.
If $N=2$ then the norm of $g'$ at $x$ is given by
$||g'(x)||_H=s$.
If $N\ge3$, then
\begin{equation}
\label{eqn:hilbertnorm}
||g'(x)||_H =
        \frac{x_N-x_{N-1}}{x_N-x_1}
        \theta\left(\frac{x_N}{x_1}\right)
        \frac{x_1^{s+1}}{E_{N-1}}
   + \frac{(x_N^s-x_1^s)x_{N-1}}{E_{N-1}}.
\end{equation}
\end{lemma}
\begin{proof}
The norm of $g'(x)$ as a map from $(\theReals^N,||\cdot||_x^H)$
to $(\theReals^N,||\cdot||_{g(x)}^H)$ is given by
\[
||g'(x)||_H = \sup_{v\in \hilbertnormball}\max_{i,k}
   \sum_{j}(h_{ij}-h_{kj})v_j,
\]
where
\[
\hilbertnormball := \Big\curly v\in\theReals^N:
            \mbox{$\max_j v_j - \min_j v_j \le 1$}\Big\ylruc.
\]
To calculate $||g'(x)||_H$ we will need to determine the sign of
$h_{ij}-h_{kj}$ for each $i,j,k\in\oneton$.
We introduce the notation
\begin{equation}
\label{eqn:elledef}
L_{ik}:=\sup_{v\in \hilbertnormball}\sum_j(h_{ij}-h_{kj})v_j.
\end{equation}

Note that $g$ is homogeneous of degree $s$, in other words
$g(\lambda x) = \lambda^s g(x)$ for all $x\in\theReals_+^N$ and $\lambda>0$.
Therefore,
\[
\sum_j x_j \frac{\partial g_i}{\partial x_j}(x)=sg_i(x)
\]
for each $i\in\curly1,\ldots,N\ylruc$.
Thus $\sum_j h_{ij}= s$ for each $i\in\curly1,\ldots,N\ylruc$, a fact
that could also have been obtained by straightforward calculation.
It follows that
\begin{equation}
\label{eqn:shift}
\sum_j(h_{ij}-h_{kj})v_j=\sum_j(h_{ij}-h_{kj})(v_j+c)
\end{equation}
for any constant $c\in\theReals$.

It is clear that an optimal choice of $v$ in~(\ref{eqn:elledef}) would be to
take $v_j:=1$ for each component $j$
such that $h_{ij}-h_{kj}>0$ and $v_j:=0$
for each component such that $h_{ij}-h_{kj}<0$.
Alternatively, we may choose $v_j:=0$ when $h_{ij}-h_{kj}>0$
and $v_j:=-1$ when $h_{ij}-h_{kj}<0$.
That the optimal value is the same in both cases follows
from~(\ref{eqn:shift}).
Also, it is easy to see that $L_{ik}=L_{ki}$.

Fix $i,k\in\curly1,\ldots,N\ylruc$ so that $i<k$.
There are four cases to consider.
\begin{itemize}
\item
{\bf Case 1.}
$1<i<k<N$.
Recall that $h_{1j}(x)=s\delta_{1j}$ and $h_{Nj}(x)=s\delta_{Nj}$.
A calculation using
Equations~(\ref{eqn:hachesstart}--\ref{eqn:hachesend}) gives
\[
E_i(x)E_k(x)(h_{i1}(x)-h_{k1}(x))
   =x_N^sx_1^{s+1}(x_k-x_i)\theta\Big(\frac{x_N}{x_1}\Big)\ge 0
\]
and
\begin{equation}
\label{eqn:norder}
E_i(x)E_k(x)(h_{iN}(x)-h_{kN}(x))
   =x_1^sx_N^{s+1}(x_k-x_i)\theta\Big(\frac{x_1}{x_N}\Big)\ge 0.
\end{equation}
We also have that
$h_{ii}(x)-h_{ki}(x)=h_{ii}(x)>0$
and $h_{ik}(x)-h_{kk}(x)= -h_{kk}(x) <0$.
So an optimal choice of $v\in \hilbertnormball$ in Equation~(\ref{eqn:elledef})
is given by $v_j:=-\delta_{jk}$.
We conclude that $L_{ik}=h_{kk}$ in this case.
\item
{\bf Case 2.}
$1=i<k<N$.
We will show that $h_{k1}(x)\le h_{11}(x)=s$.
Consider $x_1$ and $x_N$ as fixed and $x_k$ as varying in the range
$x_1\le x_k \le x_N$.
From Equation~(\ref{eqn:hachesstart}),
$h_{k1}(x) = (c_1x_k+c_2)/(c_3x_k+c_4)$,
where $c_1$, $c_2$, $c_3$, and $c_4$ depend on $x_1$ and $x_N$,
and both $c_3$ and $c_4$ are positive.
A simple calculation shows that
$c_1c_4-c_2c_3= -\theta(x_N/x_1)x_1^{s+1}x_N^s$,
which is negative.
Hence $h_{k1}$ is decreasing in $x_k$ and
takes its maximum value when $x_k=x_1$.
Here it achieves the value
\[
\frac{x_1}{x_N-x_1}\theta\Big(\frac{x_N}{x_1}\Big)
   = s-\frac{x_1^{1-s}(x_N^s-x_1^s)}{x_N-x_1}
   < s.
\]
Thus we conclude that $h_{11}(x)-h_{k1}(x)>0$. We also have that
$h_{1k}(x)-h_{kk}(x)=-h_{kk}(x) \le 0$
and $h_{1N}(x)-h_{kN}(x)=-h_{kN}(x)\ge0$.
Thus the optimal choice of $v\in \hilbertnormball$
is given by $v_j:=-\delta_{jk}$.
We conclude that in this case $L_{1k}(x)=h_{kk}(x)$.
\item
{\bf Case 3.}
$1<i<k=N$. Here $h_{i1}\ge h_{N1}=0$, $h_{ii}\ge h_{Ni}=0$, and
$h_{iN}\le h_{NN}=s$. So the optimal $v\in \hilbertnormball$ is given by
$v_j:=\delta_{j1}+\delta_{ji}$. We conclude that $L_{iN}=h_{i1}+h_{ii}$.
\item
{\bf Case 4.}
$i=1$ and $k=N$. Here $s=h_{11}\ge h_{N1}=0$ and $0=h_{1N}\le h_{NN}=s$.
Thus the optimal $v\in \hilbertnormball$ is given by
$v_j:=\delta_{1j}$. We conclude that $L_{1N}=s$.
\end{itemize}

If $N=2$ then Case 4 is the only one possible, and so $||g'(x)||_H=s$.
So, for the rest of the proof, we will assume that $N\ge3$.

We know that $h_{i1}(x)+h_{ii}(x)=s-h_{iN}(x)\ge s$ so
Case 3 dominates Case 4, that is to say $L_{iN}(x)\ge L_{1N}(x)$
for $i>1$.
Since $h_{i1}(x)\ge0$ for $i\in\curly1,\ldots,N\ylruc$, Case 3 also dominates
Cases 1 and 2, meaning that $L_{iN}(x)\ge L_{ik}(x)$ for $k<N$, $i<k$.

The final step is to maximize $L_{iN}(x)=h_{i1}(x)+h_{ii}(x)=s-h_{iN}(x)$
over $i\in\curly2,\ldots,N-1\ylruc$.
From~(\ref{eqn:norder}), $h_{mN}(x)\ge h_{nN}(x)$ for $m<n$.
Thus the maximum occurs when $i=N-1$. Recall that we have ordered 
the components of $x$ in such a way that 
$x_{N-1}$ is the second largest component of $x$.
We conclude that
\[
||g'(x)||_H = \max_{i,k:i<k}L_{ik}= h_{N-1,1}+h_{N-1,N-1}
\]
By substituting the expressions in~(\ref{eqn:hachesstart})
and~(\ref{eqn:hachesmid}),
we obtain the required formula.
\qed
\end{proof}
\begin{corollary}
\label{cor:hilbertnorm}
Let $R>0$ and $N\ge 2$.
Let $l$ be a linear functional on $\theReals^N$ such that $l(x)>0$
for all $x\in \interior \theReals_+^N$
and define $S:=\curly x\in \theReals_+^N:l(x)=1\ylruc$.
If $N=2$, then $\esssup \curly||g'(x)||_H:x\in S\ylruc=s$.
If $N\ge3$, then
\[
\esssup\curly||g'(x)||_H:\mbox{$d_H(x,\unit)\le R$, $x\in S$}\ylruc 
   = \frac{1-e^{-Rs}}{1-e^{-R}}.
\]
In both cases, the essential supremum is taken with respect
to the $N-1$-dimensional Lebesgue measure on $S$.
\end{corollary}
\begin{proof}
Note that the complement of $U\intersection S$ in $S$ has
$N-1$-dimensional Lebesgue measure zero. Using the reordering argument
in the proof of Corollary~\ref{cor:thompsonnorm}, we deduce the result
in the case when $N=2$.

The case when $N\ge3$ reduces to maximizing the
right hand side of~(\ref{eqn:hilbertnorm}) subject to the
constraints $x_1<x_{N-1}<x_N$ and $x_1/x_N\ge\exp(-R)$.
We can write the expression in~(\ref{eqn:hilbertnorm})
in the form $s+(c_1x_{N-1}+c_2)/(c_3x_{N-1}+c_4)$,
where $c_1$, $c_2$, $c_3$, and $c_4$ depend only on $x_1$ and $x_N$
and $c_1\ge0$, $c_2\le0$, $c_3\ge0$, $c_4\ge0$.
It follows that, if we view $x_1$ and $x_N$ as fixed and $x_{N-1}$ as
variable, the expression is maximized when $x_{N-1}=x_N$.
The value obtained there will be
\[
\frac{1-\left({x_1}/{x_N}\right)^s}{1-\left({x_1}/{x_N}\right)}
   =\gamma(x_1/x_N).
\]
If we recall that $\gamma$ is decreasing on $[0,1)$
and $x_1/x_N\ge\exp(-R)$, we see that
\[
||g'(x)||_H \le \frac{1-e^{-Rs}}{1-e^{-R}}.
\]
If $x_1/x_N=\exp(-R)$, then, by choosing $x\in U_{1,N}$
with $x_{N-1}$ close to $x_N$, we can
arrange that $||g'(x)||_H$ is as close as desired to this value.
\qed
\end{proof}

\begin{lemma}
\label{lem:hilbertpos}
Theorem~\ref{thm:mainhilbert} holds in the special case when
$C=\theReals_+^N$ with $N\ge 3$.
\end{lemma}
\begin{proof}
As in the proof of Lemma~\ref{lem:thompsonpos}, we may assume that
$x,y\in\interior\theReals_+^N$ and $u=\unit$.
Define $l:\theReals^N\to\theReals$ by $l(z):=\sum_{i=1}^N z_i/N$
and let $S:=\curly x\in \theReals_+^N:l(x)=1\ylruc$.
Then $l$ is a linear functional and $l(z)>0$
for all $z\in \interior \theReals_+^N$.
It is easy to check that
$\phi(s;\lambda z,\mu w)=\lambda^{1-s}\mu^s \phi(s;z,w)$
for all $\lambda,\mu>0$ and $z,w\in\interior\theReals_+^N$.
Thus
\[
d_H\Big(\phi\Big(s;\frac{u}{l(u)},\frac{x}{l(x)}\Big),
        \phi\Big(s;\frac{u}{l(u)},\frac{y}{l(y)}\Big)\Big)
   = d_H(\phi(s;u,x),\phi(s;u,y)).
\]
We also have that $d_H(x/l(x),y/l(y))=d_H(x,y)$.
Therefore we may assume that $x,y\in S$.
Let $\epsilon>0$ and define
$G:=\curly z\in S:d_H(z,\unit)< R+\epsilon\ylruc$.
It was shown in~\cite{nussbaum:hilbert} that $G$ is convex.
Also, Lemma~\ref{lem:lipschitz} states that $g$ is locally Lipschitzian.
We may therefore apply Proposition~\ref{prp:hilbert} with $f:=g$.
Since $g$ is homogeneous of degree $s$, we have that $g'(x)(x)=sg(x)$
for all $x\in G$. This, combined with the fact that $|g(x)|_{g(x)}^H=0$,
implies that $||g'(x)||_{\tilde H}=||g'(x)||_{H}$.
Using Corollary~\ref{cor:hilbertnorm}, and letting
$\epsilon$ approach zero, we deduce the required result.
\qed
\end{proof}

\begin{lemma}
\label{lem:hilbert2d}
Theorem~\ref{thm:mainhilbert} holds in the special case when
the linear span of $\curly u,x,y\ylruc$ is 1- or 2-dimensional.
\end{lemma}
\begin{proof}
If the linear span of $\curly u,x,y\ylruc$ is one-dimensional,
then all Hilbert metric distances are zero,
so assume that it is two-dimensional. 
The same argument as was used in Lemma~\ref{lem:thompson2d} shows
that it suffices to prove the result for $C=\theReals_+^2$, $u=\unit$,
and $x,y\in\interior \theReals_+^2$.
As shown in the proof of Lemma~\ref{lem:hilbertpos},
we may assume that $l(x)=l(y)=1$
where $l((z_1,z_2)):=(z_1+z_2)/2$.
We now apply Proposition 2 with $f:=g$ and
$G:=S:=\curly z\in\interior\theReals_+^2:l(z)=1\ylruc$.
Again, $||g'(x)||_{\tilde H}=||g'(x)||_{H}$ for all $x\in G$.
The result follows from the first part of Corollary~\ref{cor:hilbertnorm}.
\qed
\end{proof}

\begin{proof}[of Theorem~\ref{thm:mainhilbert}]
The proof uses Lemmas~\ref{lem:hilbertpos} and~\ref{lem:hilbert2d}
and is exactly analogous to the proof of Theorem~\ref{thm:mainthompson}.
\qed
\end{proof}

\begin{proof}[of Corollary~\ref{cor:thompson}]
We first prove the result for the case of Thompson's metric.
We will use the alternative characterization of semihyperbolicity
given in Lemma 1.2 of~\cite{semihyperbolic}.
Suppose $x,y,x',y'\in C$ are all in the same part
and are such that neither $d_T(x,x')$ nor $d_T(y,y')$ is greater than~$1$.
Let $t\in[0,\infty)$ and
write $z:=\zeta_{(x,y)}(t)$ and $w:=\phi(d_T(x,z)/d_T(x,y);x,y')$.
Observe that $d_T(y,y')\le1$ implies $|d_T(x,y)-d_T(x,y')|\le 1$.
Since $d_T(x,w)=d_T(x,y')d_T(x,z)/d_T(x,y)$, we have
\[
   |d_T(x,w)-d_T(x,z)|\le d_T(x,z)/d_T(x,y)\le 1
\]
Similar reasoning allows us to conclude that
\[
   |d_T(x,w')-d_T(x',z')|\le 1,
\]
where
$z':=\zeta_{(x',y')}(t)$ and $w':=\phi(d_T(x',z')/d_T(x',y');x,y')$.
From $d_T(x,z)=\min(t,d_T(x,y))$ and $d_T(x',z')=\min(t,d_T(x',y'))$,
we have that
\[
|d_T(x,z)- d_T(x',z')| \le |d_T(x,y)- d_T(x',y')|\le 2.
\]
So
\[
d_T(w,w')=|d_T(x,w)-d_T(x,w')|\le 4.
\]
By Theorem~\ref{thm:mainthompson}, $d_T(z,w)\le 2 d_T(y,y')\le 2$
and $d_T(z',w')\le 2 d_T(x,x')\le 2$.
The triangle inequality gives
$d_T(z,z')\le d_T(z,w)+d_T(w,w')+d_T(w',z')\le 8$.
This is the uniform bound required by the characterization
of semihyperbolicity we are using.

The proof that $C$ is semihyperbolic when endowed with Hilbert's metric 
is similar.
\qed
\end{proof}

\bibliographystyle{plain}
\bibliography{hyperbolic,topical}
 
\end{document}